\documentclass[preprint,12pt]{elsarticle}

\graphicspath{{./}}

\usepackage{amssymb, amsmath}

\usepackage[algoruled, noline, noend]{algorithm2e}

\renewcommand{\vec}{\mathbf}

\usepackage{color}

\journal{Computers \& Mathematics with Applications}

\begin{document}

\begin{frontmatter}



\title{Pattern formation in vector-valued phase fields under convex constraints}


\author[ov]{Orestis Vantzos}
\ead{ovantzos@gmail.com}
\affiliation[ov]{organization={Vantzos Research},
            city={Athens 16232},
            country={Greece}}

\begin{abstract}
In this work, a new class of vector-valued phase field models is presented, 
where the values of the phase parameters are constrained by a convex set. 
The generated phase fields feature the partition of the domain into patches of distinct phases, 
separated by thin interfaces.
The configuration and dynamics of the phases are directly dependent on the geometry and topology of
the convex constraint set, which makes it possible to engineer models of this type that exhibit desired 
interactions and patterns.

An efficient proximal gradient solver is introduced to study numerically their $L^2$-gradient flow, 
i.e.~the associated Allen-Cahn-type equation. 
Applying the solver together with various choices for the convex constraint set, 
yields numerical results that feature a number of patterns observed in nature and engineering, 
such as multiphase grains in metal alloys, traveling waves in reaction-diffusion systems, and vortices in magnetic materials.

\end{abstract}



\begin{keyword}
phase field models\sep obstacle potential\sep convex optimization\sep proximal gradient method\sep pattern formation

\MSC[2010] 00-01\sep  99-00

\end{keyword}

\end{frontmatter}




\section{Introduction}
\label{sec:intro}

\paragraph{Overview} 

In the rest of section \ref{sec:intro}, the celebrated Ginzburg-Landau functional $\int_\Omega \frac{\epsilon}{2}\lvert \nabla u\rvert^2 + \frac{1}{\epsilon}(1-u^2)^2\,\mathrm{d}x$,
introduced in material science \cite{cahn1958free} in the form of the Allen-Cahn and Cahn-Hilliard equations
for the study of two-phase materials, is presented. 
Its extensions into vector-valued phase fields $\vec u$, and obstacle potentials which restrict the possible values of the phase field,
are also presented. These extensions are meant to facilitate the modeling of materials with more than two phases.

In section \ref{sec:theory}, a novel modification of the functional is introduced, 
characterized by the confinement of the values of a vector-valued phase field $\vec{u}$ within a general convex constraint set $C$.
The optimality conditions of this new functional are presented, 
and used in turn to derive the properties of the phases and interfaces of the phase field $\vec{u}$
as a function of the geometry and topology of the constraint set $C$. 

The development of a robust numerical scheme for the $L^2$-gradient flow of the functional (the corresponding `Allen-Cahn' evolution equation)
is presented next in section\ref{sec:solver},
based on ideas from the theory of gradient flows in metric spaces \cite{ambrosio2008gradient} and proximal minimization for monotone operators \cite{bauschke2011convex}.
Using the scheme to evolve the phase field $\vec u$ from initial random conditions and under a range of different choices for the constraint set $C$ (see Fig.~\ref{fig:obstacles}), 
yields a series of numerical results (Fig.~\ref{fig:twophases}--Fig.~\ref{fig:twave}) in section \ref{sec:num_results}.

\begin{figure}
 \centering
 \frame{\includegraphics[height=140pt]{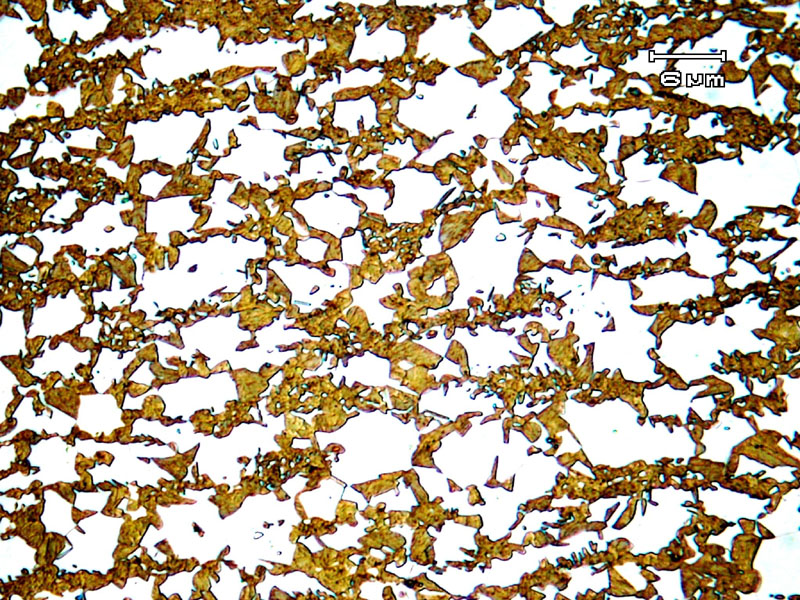}}
 \frame{\includegraphics[height=140pt]{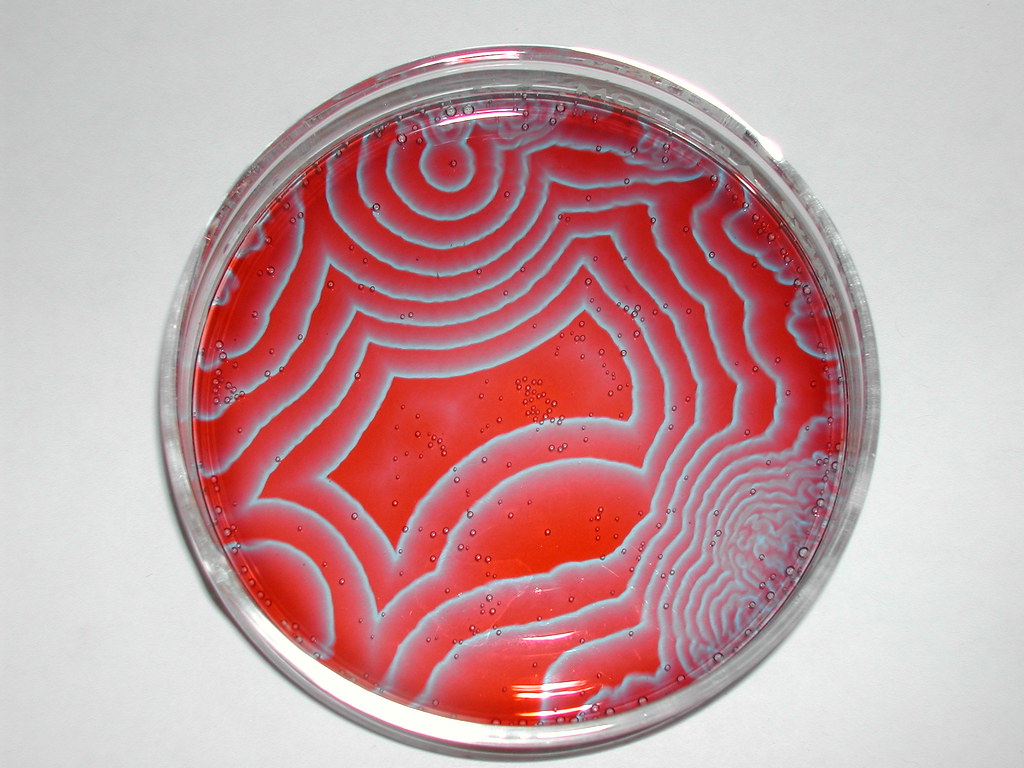}}\\[3pt]
 \frame{\includegraphics[height=151pt]{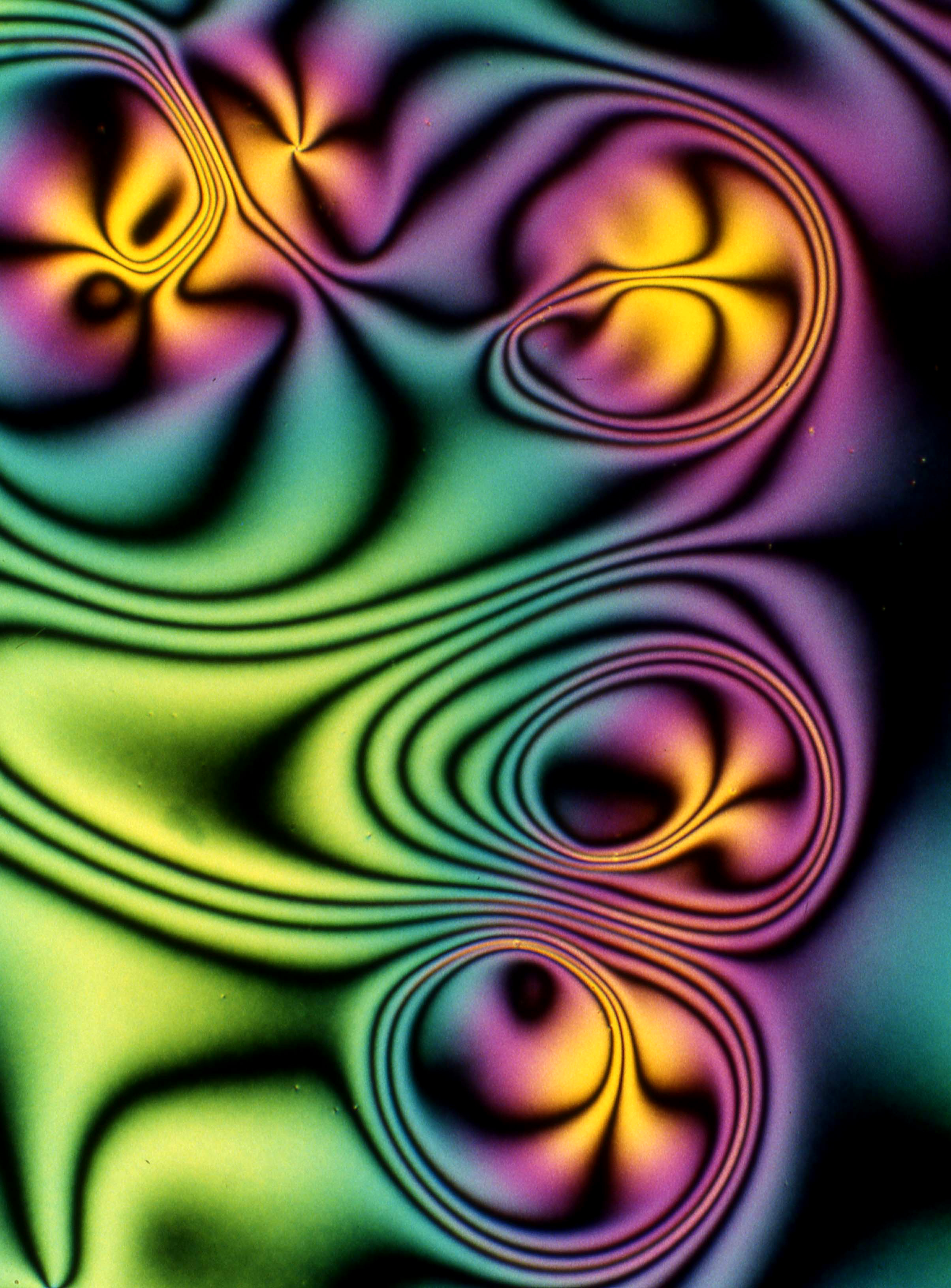}}
 \frame{\includegraphics[height=151pt]{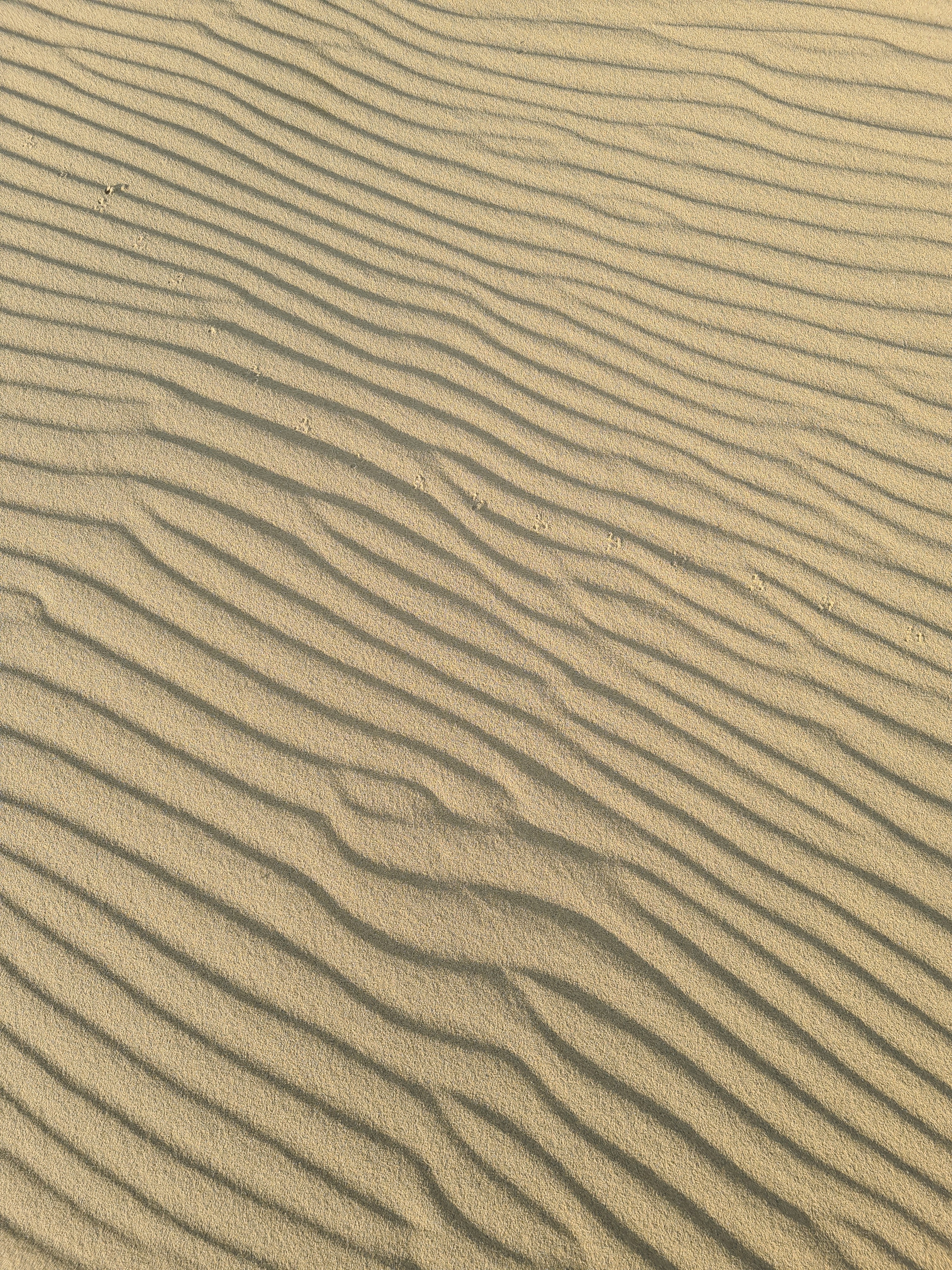}}
 \frame{\includegraphics[height=151pt]{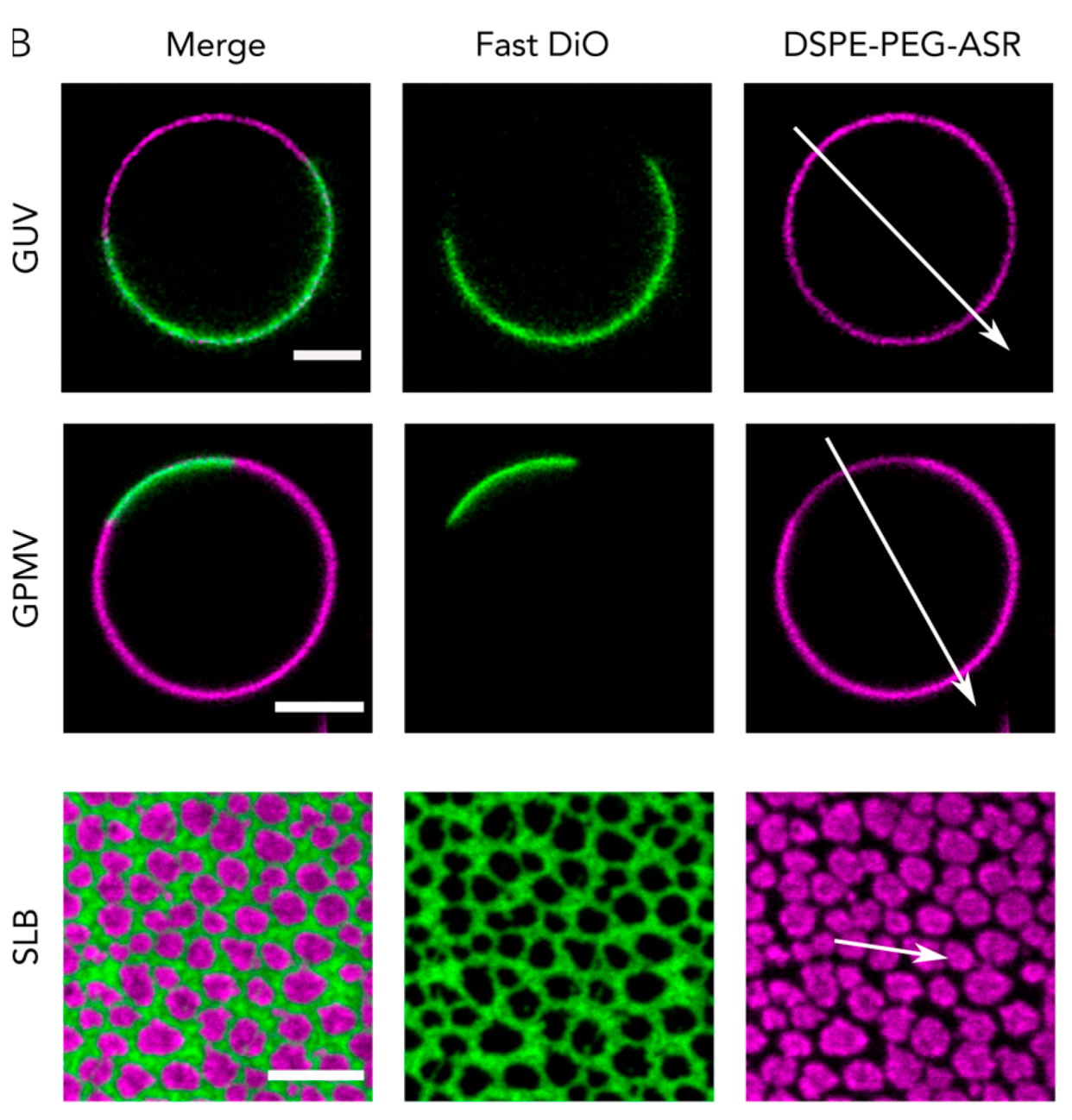}}\\
 \caption{\emph{Patterns in Nature and Engineering}. From top left, clockwise: 
 i) Grains in a dual-phase steel alloy (Credit: Dr Amar K De, courtesy of University of Cambridge).
 ii) Traveling waves during a Belousov-Zhabotinsky reaction (Credit: Stephen Moris, courtesy of Flickr).
 iii) Phase-separation in GUVs, a model for cell membranes (Credit: Taras Sych et al., courtesy of MDPI).
 iv) Wave patterns in desert sand (Credit: Olga~K, courtesy of Pexels).
 v) Topological point defects in a thin liquid crystal film (Credit: Oleg Lavrentovich, courtesy of NSF).
 }
 \label{fig:examples}
 \end{figure}

The numerical results reveal the remarkable ability of this generalized functional 
to model a range of patterns far wider than the original context of grain formation in multiphased materials,
such as those presented in Fig.~\ref{fig:examples}. 
In addition to a discussion of the immediate applications in terms of modeling the processes that generate these patterns, 
especially when coupled with equations such as the Navier-Stokes, 
one can find in section \ref{sec:future} a selection of ideas for future relevant work, in Analysis, Numerics and Computer Science.

\paragraph{The Phase Field Method}

\begin{figure}
  \centering
 \includegraphics[width=0.8\textwidth]{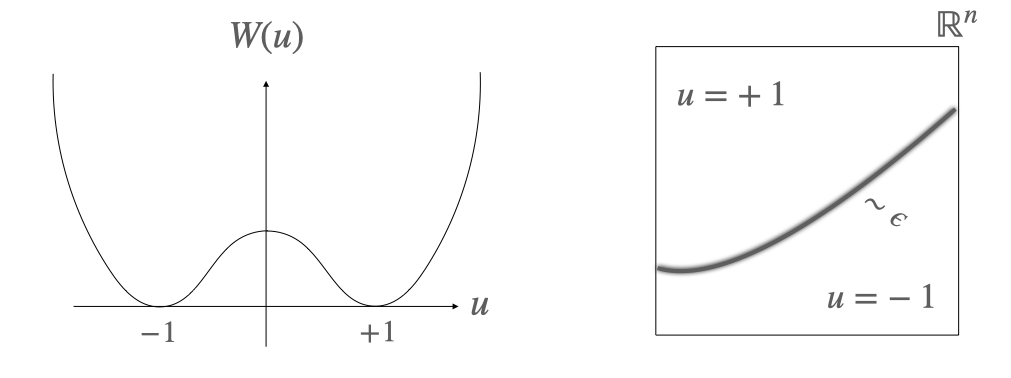}
 \caption{Double-well potential $W$ (left). 
 Partition of domain into phases separated by interfaces of thickness $O(\epsilon)$, typical for solutions of the classic Allen-Cahn equation (right).}
 \label{fig:allencahn}
 \end{figure}

Consider the \emph{Allen-Cahn} equation
\begin{equation}
 \frac{\partial u}{\partial t} = - \epsilon \Delta u + \epsilon^{-1}\nabla W(u),\quad W(u)=\frac{1}{4}(1-u^2)^2
\end{equation}
where $u:\Omega\times [0,T]\rightarrow \mathbb{R}$, $\Omega\subset \mathbb{R}^n$ and $T>0$, and $\epsilon$ is a small positive constant.
This is a well-studied 2nd-order evolution PDE, 
whose solutions feature (for large enough $t$) partitions of the domain into regions with constant $u=\pm 1$, 
representing two distinct \emph{phases}, separated by interfaces with thickness $\sim \epsilon$ \cite{chen1994spectrum}. 
See Fig.~\ref{fig:allencahn} for a sketch of this idea, 
and the top left image in Fig.~\ref{fig:examples} for an example of the real world process that motivated Allen and Cahn.

One can derive this equation as the \emph{$L^2$-gradient flow} of the \emph{Ginzburg-Landau functional}
\begin{equation}
 \mathcal{G}(u) = \int_\Omega \frac{\epsilon}{2}\lvert \nabla u\rvert^2 + \frac{1}{\epsilon}W(u)\,\mathrm{d}x\,,
\end{equation}
introduced in \cite{cahn1958free} in the context of material science. 
This leads to the classic result that the interface between the phases follows, in a suitable sense, motion by mean-curvature \cite{ilmanen1993convergence}.

The idea to approximate a partition of a domain into separate phases via a smooth function, 
known as the \emph{phase field method} \cite{emmerich2003diffuse}, has found many practical applications, 
where the Allen-Cahn, or the related 4th-order Cahn-Hilliard PDE, is coupled with other equations 
such as the Navier-Stokes \cite{giorgini2019uniqueness} or elasticity equations \cite{penzler2012phase}.

\paragraph{Vector-Valued Phase Fields}

The need to extend the phase field method to cases where there are more than two phases led to the study of \emph{vector-valued} equations of Allen-Cahn type
\begin{equation}\label{eq:AllenCahn}
 \frac{\partial u_i}{\partial t} = - \epsilon \Delta u_i + \epsilon^{-1}\frac{\partial W(\vec u)}{\partial u_i}\end{equation}
 where $\vec{u}:\Omega\times[0,T]\rightarrow \mathbb{R}^m$ is the vector-valued phase field, 
 and $W:\mathbb{R}^m\rightarrow \mathbb{R}$ is a \emph{multi-well potential} with minima $\vec{U}_1,\ldots,\vec{U}_k \in\mathbb{R}^m$ representing $k$ distinct phases.
 The dynamics of the vector-valued Allen-Cahn are sensitive to the exact shape of the potential $W$, 
 for instance the existence and location of \emph{heteroclinics} that connect the various minima \cite{alikakos2008connection}, \cite{bates2017multiphase}.
 
\paragraph{Obstacle Potentials} 

\begin{figure}
  \centering
 \includegraphics[width=0.9\textwidth]{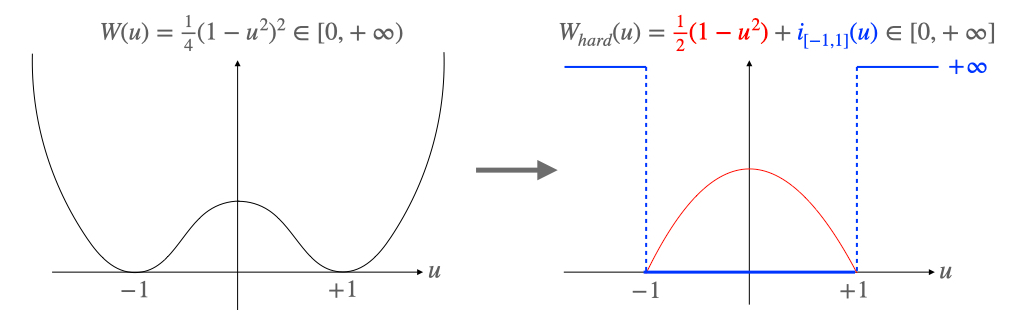}
 \caption{The idea behind \emph{obstacle potentials} is to replace a higher-order polynomial multi-well potential (left),
 with the sum of a concave quadratic polynomial plus the characteristic function of a suitable convex set (right).
 The minima (wells) of the soft potential correspond to extremal points of the convex constraint set.}
 \label{fig:hardpotential}
 \end{figure}

A different approach to multiphase dynamics is to constrain the values of the phase field to the \emph{$k$-dimensional Gibbs simplex $D^k$}:
\begin{equation}
 \vec{u}\in D^k \text{ iff } \sum_{i=1}^k u_i = 1 \text{ and } u_i\geq 0\,.
\end{equation}
The vertices $\vec{V}_1,\ldots,\vec{V}_k$ of the simplex, 
where a single component $u_i$ takes the value 1 and all the rest vanish, represent then $k$ distinct phases \cite{emmerich2003diffuse}. 
Contrary to the unconstrained case, where the multi-well potentials $W$ are by necessity complicated functions, e.g.~high-degree multivariate polynomials, 
the concave quadratic potential $W(\vec u)  = \sum_{i=1}^k u_i(1-u_i) $, constrained to the simplex $D^k$, 
achieves its minima exactly at the $k$ vertices of the simplex. 
In Fig.~\ref{fig:hardpotential}, one can see how this process of `hardening' the potential works for a 1-dimensional simplex, the interval $[-1,1]$.

\section{Ginzburg-Landau-type Functionals with Convex-Obstacle Potential}
\label{sec:theory}

\paragraph{The Generalized Functional} 

We generalize the concepts of the previous section by introducing a \emph{generalized Ginzburg-Landau functional} $\mathcal{G}_C$ 
for a vector-valued phase field $\vec{u}:\Omega\subset \mathbb{R}^n\rightarrow\mathbb{R}^m$  
with values restricted to a (bounded and closed) convex set $C\subset\mathbb{R}^m$:
\begin{equation}\label{eq:GL}
 \mathcal{G}_C(\vec u) := \int_\Omega \frac{\epsilon}{2}\sum_{i=1}^m \lvert \nabla u_i\rvert^2 + \frac{1}{\epsilon} W(\vec u) + i_C(\vec u)\,\mathrm{d}x^n\,.
\end{equation}
The potential $W:\mathbb{R}^m\rightarrow \mathbb{R}$ is assumed to be \emph{differentiable, concave and non-negative in $C$}, such as the quadratic potential
\begin{align}
 &W(\vec y)  = \frac{1}{2}(d_C^2 - \lvert \vec y\rvert^2)\\
 &d_C = \max_{\vec y\in C}\,\lvert \vec y\rvert\,.
\end{align} 
The extended-valued function
\begin{equation}
 i_C(\vec y) = \begin{cases}
  0, & \vec y\in C\\
  \infty, & \vec y \notin C
  \end{cases}
\end{equation}
is \emph{the indicator function of the constraint set $C$}.
One can directly check that the functional is non-negative, $\mathcal{G}_C(\vec u)\geq 0$, 
and that $\mathcal{G}_C(\vec u)<\infty$ iff $\vec{u}(\vec x)\in C$ almost everywhere in $\Omega$.

\begin{figure}
 \centering
\includegraphics[width=0.9\textwidth]{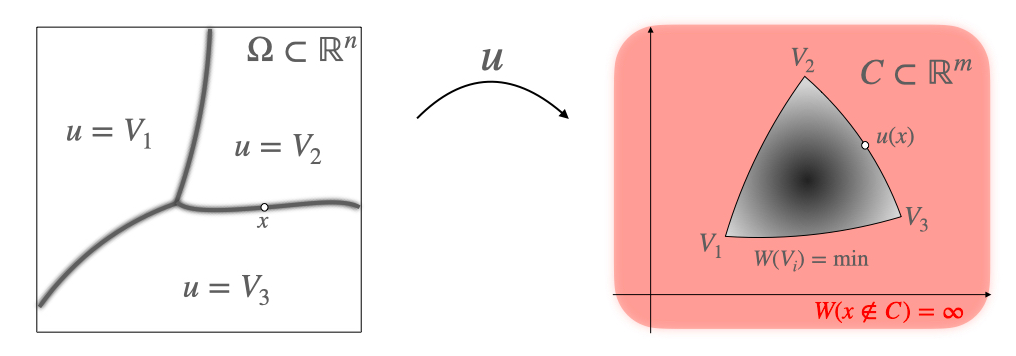}
\caption{Partition of the domain into multiple phases, typical of the minimizers of the generalized Ginzburg-Landau functional (left), 
and the corresponding convex constraint set (right).}
\label{fig:genGL}
\end{figure}

\paragraph{Optimality Conditions} 

To derive the first-order optimality conditions for the functional $\mathcal{G}_C$, 
we need to consider small perturbations around a feasible phase field $\vec u$. 
Since the perturbed field needs to also be feasible, 
we consider perturbations of the form $\vec u\rightarrow \vec u_\theta = \vec u + \theta(\vec w-\vec u)$ with $\vec w$ feasible 
and $0\leq \theta \ll 1$, i.e.~interpolations in the direction of another feasible field. 
The variation of the functional under such a perturbation is then
\begin{multline}
 \mathcal{G}_C(\vec u_\theta) =  \mathcal{G}_C(\vec u)\\
  + \theta \int_\Omega\left\{ \epsilon \sum_{i=1}^m\nabla u_i\cdot \nabla (w_i-u_i) + \frac{1}{\epsilon} \sum_{i=1}^m \frac{\partial W(\vec u)}{\partial u_i}(w_i-u_i)\right\}\mathrm{d}x^n + \operatorname{O}(\theta^2)
\end{multline}
The first order optimality condition, i.e.~that $\mathcal{G}_C(\vec u_\theta)\geq\mathcal{G}_C(\vec u) $ to first order in $\theta$, 
can be written then (using the vector $L^2$-product notation $(\vec u,\vec v)=\int_\Omega \vec u\cdot\vec v\,\mathrm{d}x^n$) as the following \emph{variational inequality}:
\begin{equation}\label{eq:optWeak}
 \epsilon\sum_{i=1}^m\left(\nabla u_i,\nabla(w_i-u_i)\right) - \frac{1}{\epsilon}\left(\nabla W(\vec u),\vec w-\vec u\right) \geq 0, \quad\forall \vec w:\vec w(x)\in C\text{ a.e.~in }\Omega
\end{equation}
Alternatively, we can rewrite the optimality condition in strong form as a \emph{differential inclusion}:
 \begin{equation}\label{eq:optStrong}
 \epsilon\Delta\vec u(x) - \frac{1}{\epsilon}\nabla W(\vec u(x))\in \partial i_C(\vec u(x))\text{ a.e.~in }\Omega
\end{equation}
where the Laplacian operator is applied component-wise $(\Delta \vec u)_i = \Delta u_i$. 
The \emph{subdifferential} $\partial f(\vec x)$ of a convex extended-valued function $f:\mathbb{R}^n\rightarrow \mathbb{R}\cup\infty$ at a point $\vec x\in \mathbb{R}^n$ 
is defined as the set of $\vec v\in\mathbb{R}^n$ (the \emph{subgradients}) such that $f(\vec y)-f(\vec x)\geq \vec v\cdot(\vec y-\vec x)$ for any $\vec y\in  \mathbb{R}^n$. 
In the particular case where the function $f$ is the indicator function $i_C$ of the convex set $C$, 
it can be shown that for any $\vec x\in C$, $\vec v\in\partial i_C(\vec x)$ iff the variational inequality $\vec v\cdot(\vec y-\vec x)\leq 0$, for any $\vec y\in C$, is satisfied.
Combining this with integration by parts connects \eqref{eq:optStrong} and \eqref{eq:optWeak}. 
From a geometrical point of view, we can identify the subdifferential $\partial i_C(\vec u)$ with the \emph{normal cone} $N(C;\vec u)$ of the convex set $C$ at the point $\vec u\in C$. 
In the case where the boundary $\partial C$ admits a continuous (outward-pointing) normal field $\vec n$ in the neighborhood of $\vec u$, 
then the normal cone is simply $N(C;\vec u) = \{\lambda \vec n(\vec u)\}_{\lambda\geq 0}$.

\paragraph{Phases}  In line with the classic theory of \emph{phase separation}, as modelled by the Allen-Cahn equation \eqref{eq:AllenCahn}, 
the functional $\mathcal{G}_C$ \eqref{eq:GL} favors (for small values of the parameter $\epsilon$) phase fields $\vec u:\Omega\rightarrow C$ 
that correspond to partitions of the domain $\Omega$ into areas where the potential $W$ is minimal (over $C$) and the spatial variation of $\vec u$ is small. 
This allows us to ignore the Laplacian term in \eqref{eq:optStrong}, and the optimality condition reduces to 
\begin{equation}
 -\nabla W(\vec u) \in N(C;\vec u) \Leftrightarrow \nabla W(\vec u)\cdot(\vec w-\vec u)\geq 0\,,
\end{equation}
i.e.~the potential at a minimal point increases in any feasible direction that stays inside the constraint set.
Since the potential $W$ is assumed concave, it can only have local maxima in the interior of the convex set $C$. 
All the minima of $W$ over $C$ lie on the boundary $\partial C$ then. 

Isolated minima $\vec U_1,\vec U_2,\ldots$ correspond to uniform phases where $\vec u(x)=\vec U_k=const$, 
and so $\Delta\vec u = 0$, equivalent to the classic Ginzburg-Landau functional's two phases where $u(x)=\pm 1=const$. 
Contrary to the Ginzburg-Landau, we might also have a connected set $S\subset\partial C$ of minima, such that $W(\vec u)  = const$ for any $\vec u\in S$. 
We can characterize such a phase as a \emph{harmonic map} from an open subset of $\Omega$ to $S$ (considered as a submanifold of $\mathbb{R}^m$). 
Indeed, since $\vec u(x)\in S$ and so $W(\vec u(x))=const$ everywhere, 
the phase field $\vec u$ is a critical point of the functional $\mathcal{G}_C$ only if it is a critical point of the functional 
\begin{multline}
 \int_\Omega \frac{1}{2}\sum_{i=1}^m \lvert \nabla u_i\rvert^2 \,\mathrm{d}x^n = \int_\Omega \frac{1}{2}\sum_{j=1}^n \lvert \frac{\partial \vec u}{\partial x_j}\rvert^2 \,\mathrm{d}x^n\\
 = \int_\Omega \frac{1}{2}\sum_{j=1}^n \lvert \pi_{T_{\vec u}S}\frac{\partial \vec u}{\partial x_j}\rvert^2 \,\mathrm{d}x^n = \int_\Omega \frac{1}{2} \lvert \mathrm{d}(id_{S}^{-1}\circ\vec u)\rvert^2\,\mathrm{d}x^n \,,
\end{multline}
where $\pi_{T_{\vec u}S}$ is the orthogonal projection from $\mathbb{R}^m$ to the tangent space of $S$ at $\vec u\in S$, and $id_S$ is the trivial embedding of $S$ into $\mathbb{R}^m$. 
A corresponding phase then can be described by the following conditions: $\vec u(x)\in S$ almost everywhere, together with the system of non-linear second-order elliptic PDEs 
\begin{equation}
 -\sum_{j=1}^n \pi_{T_{\vec u} S}\frac{\partial}{\partial x_j} \left(\pi_{T_{\vec u} S} \frac{\partial \vec u}{\partial x_j}\right)=0 \,.
\end{equation}
In particular, if $C\subset\mathbb{R}^m$ and $S$ is an (m-1)-dimensional hypersurface with unit normal $\vec n(\vec u)$, 
then $\pi_{T_{\vec u}S} =\operatorname{id} - \frac{\vec{n}(\vec u)\vec{n}(\vec u)^T}{\lvert \vec n(\vec u)\rvert^2}$, 
whereas if $S$ is a one dimensional curve with unit tangent vector $\vec t(\vec u)$ then $\pi_{T_{\vec u}S} = \vec t(\vec u) \vec t(\vec u)^T$.

\paragraph{Interfaces} As a first approximation, we consider a planar interface between two uniform phases, $\vec u(x)=\vec u_1$ over the half-space $x_1\leq 0$ 
and $\vec u(x)=\vec u_2$ over the half-space $x_1\geq \epsilon \Lambda>0$, that correspond to two isolated minima of $W(\vec u)$ over $C$. 
The profile of the interface then is a function $\vec U:[0,\Lambda]\rightarrow C$, $\vec U(0)=\vec u_1$ and $\vec U(T)=\vec u_2$, 
such that $\vec u(x)=\vec U(\frac{x_1}{\epsilon})$ when $0\leq x_1\leq \epsilon \Lambda$. 
We assume moreover that $W(\vec U(\xi))>0$ for any intermediate point $\xi\in(0,\Lambda)$, 
i.e.~the path connecting $\vec u_1$ and $\vec u_2$ does not go through any other minimal point. 

Integrating the generalized Ginzburg-Landau functional \eqref{eq:GL} over a cylindrical cross-section of the interface of the form $\Omega_\epsilon = [0,\Lambda]\times D$, 
where $D$ is a subset of the $x_1=0$ plane, we have
\begin{multline}
 \int_{\Omega_\epsilon} \frac{\epsilon}{2}\sum_{i=1}^m\lvert\nabla u_i\rvert^2 + \frac{1}{\epsilon}W(\vec u) + i_C(\vec u)\,\mathrm{d}x^n\\
= \int_D \int_0^\Lambda \frac{1}{2}\lvert \vec U'(\xi)\rvert^2 + W(\vec U(\xi))\,\mathrm{d}\xi\,\mathrm{d}x^{n-1}
\end{multline}

We consider an arbitrary feasible perturbation $\vec U\rightarrow \vec U_\theta =  \vec U + \theta(\vec w-\vec U)$, with $0<\theta\ll 1$, $\vec w(\xi)\in C$ and $\vec w(0)=\vec w(T)=0$:
\begin{multline}
 \int_0^\Lambda \frac{1}{2}\lvert \vec U_\theta'(\xi)\rvert^2 + W(\vec U_\theta(\xi))\,\mathrm{d}\xi\\
 = \theta \int_0^\Lambda (-\vec U''(\xi) + \nabla W(\vec U(\xi)))\cdot(\vec w(\xi)-\vec U(\xi))\,\mathrm{d}\xi + \operatorname{O}(\theta^2)
\end{multline}
The first order optimality condition for the profile $\vec U$, under the convex constraint, can be written then as
\begin{multline}\label{eq:optU}
 (-\vec U''(\xi) + \nabla W(\vec U(\xi)))\cdot(\vec w-\vec U(\xi))\geq 0,\,\forall \vec w\in C\\
 \Rightarrow \vec U''(\xi) - \nabla W(\vec U(\xi)) \in N(C;\vec U(\xi))
\end{multline}
for (almost) all $\xi\in[0,\Lambda]$.

To proceed, we assume furthermore that the profile $\vec U$ corresponds to a path $\gamma:[0,L]\rightarrow C$ 
joining $\gamma(0)=\vec u_1$ to $\gamma(L)=\vec u_2$ through the constraint set $C$, parametrized by its arc-length $s$, 
via the reparametrization $\vec U(\xi) = \gamma(s(\xi))$ with $s'(\xi)>0$ a.e.~in $(0,\Lambda)$. 
Via a change of variables from $\xi$ to $s$, the interface energy can be written as
\begin{equation}\label{eq:optE}
 E = \int_0^\Lambda \frac{1}{2}\lvert \vec U'(\xi)\rvert^2 + W(\vec U(\xi))\,\mathrm{d}\xi
= \int_0^L \left(\frac{1}{2\xi'(s)} + \xi'(s)W(\gamma(s))\right)\,\mathrm{d}s
\end{equation}
which is minimized when 
\begin{equation}
 \xi'(s)=\lvert 2W(\gamma(s))\rvert^{-1/2} \Rightarrow s'(\xi) = \sqrt{2W(\vec U(\xi))}
\end{equation}
From this we can calculate both the interfacial energy $E$ and its (rescaled) thickness $\Lambda$:
\begin{align}
 &E = \int_0^L \sqrt{2W(\gamma(s))}\,\mathrm{d}s\\
 \label{eq:optLambda}&\Lambda = \int_0^L \xi'(s)\,\mathrm{d}s = \int_0^L \lvert 2W(\gamma(s))\rvert^{-1/2}\,\mathrm{d}s
\end{align}
Going back to the first order optimality condition \eqref{eq:optU}, it can be rewritten in terms of $\gamma$ as
\begin{equation}\label{eq:optG}
 \sqrt{2 W(\gamma(s))}\gamma''(s) - (\operatorname{id}-\vec t(s)\vec t(s)^T)\nabla W(\gamma(s)) \in N(C;\gamma(s))
\end{equation}
where $\vec t(s) = \gamma'(s)$ is the unit tangent vector of $\gamma$.

We consider finally the special case where the path $\gamma$ lies entirely on the boundary $\partial C$ of the constraint set 
and, moreover, there exists a continuous unit normal $\vec n$ such that  $N(C;\gamma(s)) = \{\lambda\vec n(s)\}_{\lambda\geq 0}$ for any $s\in(0,L)$. 
There exists then a unit tangent vector $\vec b(s)$, such that $\nabla W(\gamma(s)) \in \operatorname{span}(\{\vec t(s), \vec n(s),\vec b(s)\})$. 
The triplet $\{\vec t(s), \vec n(s),\vec b(s)\}$ serves as a \emph{Darboux frame} \cite{do2016differential} for $\gamma$, 
and the acceleration can be decomposed as $\gamma''(s)=\kappa_n(s)\vec n(s) + \kappa_g \vec b(s)$, where $\kappa_n, \kappa_g$ is the normal and geodesic curvature respectively. 
This allows us to reduce the optimality condition \eqref{eq:optG} to $\gamma''(s)\in\operatorname{span}\{\vec n(s),\vec b(s)\}$ and 
\begin{align}
 & \sqrt{2W(\gamma(s))}\,\kappa_n - \nabla W(\gamma(s))\cdot \vec n(s) \geq 0\,,\\
 & \sqrt{2W(\gamma(s))}\,\kappa_g - \nabla W(\gamma(s))\cdot \vec b(s) = 0\,.
\end{align}
Taking into account that the normal curvature $\kappa_n\leq 0$ since $C$ is convex, 
and therefore the normal derivative $\frac{\partial W}{\partial n} = \nabla W\cdot \vec n$ should also be non-positive over an optimal path, 
we can rewrite the first optimality condition above as
\begin{equation}
 \lvert \kappa_n \rvert \leq \left\lvert \frac{\partial \sqrt{2W}}{\partial n}\right\rvert\,.
\end{equation}
This yields the following important insight: 
\emph{the interface $\vec U(\xi)$ between two phases $\vec u_1,\vec u_2\in \partial C$ only takes values from the boundary $\partial C$ of the constraint set, 
as long as the gradient of the potential $W$ is large enough compared to the curvature of the boundary.}

\section{Proximal Gradient Solver for the $L^2$-Gradient Flow}
\label{sec:solver}

\paragraph{Minimizing Movements} 
We are interested in studying numerically the Allen-Cahn analogue for the generalized Ginzburg-Landau functional \eqref{eq:GL}, i.e.~the $L^2$-gradient flow of the functional $\mathcal{G}_C$. 
Motivated by the method of \emph{minimizing movements} for gradient flows \cite{ambrosio2008gradient}, 
we consider the following variational time-discrete scheme, which yields an updated phase field $\vec u^{k+1}$ at time $t=t^k+\tau$, given the phase field $\vec u^k$ at time $t=t^k$:
\begin{equation}\label{eq:minmov}\begin{aligned}
 \vec u^{k+1} &= \operatorname*{argmin}_{\vec u\in \vec u^0 + H_0^1(\Omega)^m}\left\{\frac{1}{2\tau}\lVert \vec u-\vec u^k\rVert_{L^2}^2 + \mathcal{G}_C(\vec u) \right\}\\
 &=  \operatorname*{argmin}_{\vec u\in \vec u^0 + H_0^1(\Omega)^m}\left\{\int_\Omega \frac{1}{2\tau}\lvert \vec u-\vec u^k\rvert^2 \right.\\ 
 &\hspace{130pt}\left.+\frac{\epsilon}{2}\sum_{i=1}^m\lvert\nabla u_i\rvert^2 + \frac{1}{\epsilon}W(\vec u) + i_C(\vec u)\,\mathrm{d}x^n\right\}
\end{aligned}\end{equation}
where $\vec u^0$ is the initial phase field at $t=0$. The phase field $\vec u^0$ also encodes the Dirichlet boundary conditions at the boundary of the domain $\Omega$.

Assuming that there exists a time step $\tau>0$ small enough so that $\frac{1}{2\tau}\lvert \vec u-\vec u^k\rvert^2 + \frac{1}{\epsilon}W(\vec u)$ is convex in $C$, 
it is straightforward to show that the objective functional of \eqref{eq:minmov} is proper, convex and lower semi-continuous in $H^1(\Omega)^m$. 
Noting that any minimizer of the functional \eqref{eq:minmov} satisfies the inequality
\begin{equation}
  \frac{1}{2\tau}\lVert \vec u-\vec u^k\rVert_{L^2}^2 + \mathcal{G}_C(\vec u) \leq \mathcal{G}_C(\vec u^k)\,,
\end{equation}
we can show that a minimizing sequence is essentially bounded in $H^1(\Omega)^m$, 
since the inequality above implies the bound
\begin{equation}
\frac{\epsilon}{2} \sum_{i=1}^m \lVert \nabla u_i\rVert^2 \leq \mathcal{G}_C(\vec u) \leq \mathcal{G}_C(\vec u^k).
\end{equation} 
Combined with the lower bound $\frac{1}{2\tau}\lVert \vec u-\vec u^k\rVert_{L^2}^2 + \mathcal{G}_C(\vec u) \geq 0$, 
we deduce the existence of $\vec u^{k+1}$ as defined in \eqref{eq:minmov}, see \cite{dacorogna2007direct}. 
If moreover we assume that $\frac{1}{2\tau}\lvert \vec u-\vec u^k\rvert^2 + \frac{1}{\epsilon}W(\vec u)$ is strongly convex, for small enough $\tau$, then the minimizer $\vec u^{k+1}$ is unique.

\paragraph{A Splitting Algorithm} To derive a practical algorithm for the solution of the optimization problem \eqref{eq:minmov}, 
we will make use of a variant of the proximal gradient method, a splitting algorithm introduced by Tseng in \cite{tseng2000modified} (see Prop.~27.13 in \cite{bauschke2011convex}). 
Assuming that $\mathcal{F}$ and $\mathcal{G}$ are proper, convex and lower semi-continuous, with $\mathcal{G}'$ $1/\beta$-Lipschitz continuous, 
and that $C$ is a closed convex set, then Tseng's algorithm is the following iteration: starting with $u^{(0)} = \vec u^k$, repeat
\begin{subequations}\label{eq:tseng}\begin{align}
   &\vec x^{(m)} = \vec u^{(m)} - \gamma \mathcal{G}'(\vec u^{(m)})\\
   &\vec y^{(m)} = \operatorname{prox}_{\gamma \mathcal{F}} \vec x^{(m)}\\
   &\vec z^{(m)} = \vec y^{(m)} - \gamma \mathcal{G}'(\vec y^{(m)})\\
   &\vec u^{(m+1)} = \operatorname{proj}_C(\vec u^{(m)} - \vec x^{(m)} + \vec z^{(m)})
\end{align}\end{subequations} 
with a fixed $\gamma\in(0,\beta)$. At the limit, $\lVert\vec y^{(m)}-\vec u^{(m)}\rVert \rightarrow 0$, 
and $\vec y^{(m)}$, $\vec u^{(m)}$ converge weakly to a minimizer of $\mathcal{F}+\mathcal{G}$ with values in $C$.
Under more strict conditions on the convexity of $\mathcal{F}$ or $\mathcal{G}$ (see again Prop.~27.13 in \cite{bauschke2011convex}),
there is in fact a unique feasible minimizer of $\mathcal{F}+\mathcal{G}$ and $\vec y^{(m)}$, $\vec u^{(m)}$ converge strongly to it.

We split the functional in $\mathcal{F}(\vec u)=\frac{\epsilon}{2}\sum_i\lVert\nabla u_i\rVert^2$ and $\mathcal{G}(\vec u)=\frac{1}{2\tau}\lVert\vec u-\vec u^k\rVert^2 + \frac{1}{\epsilon}\int_\Omega W(\vec u)\,\mathrm{d}x 
\Rightarrow \mathcal{G}'(\vec u) = \frac{1}{\tau}(\vec u-\vec u^k) + \frac{1}{\epsilon}\nabla W(\vec u)$. 
For the Lipschitz continuity, we have $\mathcal{G}'(\vec u)- \mathcal{G}'(\vec v) = \frac{1}{\tau}(\vec u-\vec v) + \frac{1}{\epsilon}(\nabla W(\vec u)-\nabla W(\vec v))$, 
and so if $\nabla W$ is $1/\omega$-Lipschitz continuous over $C$, then $\mathcal{G}'$ is $1/\beta$-Lipschitz continuous with $\beta=(\tau^{-1}+(\epsilon\omega)^{-1})^{-1}\leq\min(\tau,\epsilon\omega)$. 
Furthermore, we can show that the $1/\omega$-Lipschitz-continuity of $\nabla W$ implies that $\mathcal{G}$ is \emph{strongly convex}, for $\tau$ small enough, 
and therefore the $\vec u^{(m)}$ converge strongly to the unique minimizer of the problem. 
Indeed, if $\lvert \nabla W(\vec u)-\nabla W(\vec v)\rvert \leq \omega^{-1}\lvert \vec u-\vec v\rvert$, then
\begin{multline}
  (\nabla W(\vec u)-\nabla W(\vec v))\cdot(\vec u-\vec v) \geq -\lvert\nabla W(\vec u)-\nabla W(\vec v)\rvert\lvert \vec u - \vec v\rvert\\
  \geq -\frac{1}{\omega}\lvert \vec u-\vec v\rvert^2
\end{multline}
which implies that $W(u) + \frac{1}{2\omega}\lvert \vec u\rvert^2$ is convex. 
From this, we can deduce that $\mathcal{G}$ is strongly convex when $\frac{1}{\tau}>\frac{1}{\omega\epsilon}\Rightarrow \tau < \omega\epsilon$.
Combining all the conditions over the various parameters, we conclude that a sufficient condition for the iteration to converge is $0<\gamma<\tau<\epsilon\omega$. 

For the proximal operator of $\mathcal{F}$ in $L^2$, we have by definition
\begin{multline}
  \vec u^*:=\operatorname{prox}_{\gamma\mathcal{F}}(\vec u) = \operatorname*{argmin}_{\vec v\in L^2(\Omega)^m} \left\{\frac{1}{2\gamma}\lVert \vec u-\vec v\rVert^2 + \mathcal{F}(\vec v) \right\} \Rightarrow\\ 
  \langle \frac{1}{\gamma}(\vec u^*- \vec u), \vec v \rangle + \sum_{i=1}^m\langle \epsilon \nabla u_i^*,\nabla v_i\rangle = 0,\quad\,\vec v\in H_0^1(\Omega)^m\\
  \Rightarrow \operatorname{prox}_{\gamma\mathcal{F}}(\vec u) = (id-\epsilon\gamma\Delta)^{-1}\vec u
\end{multline}
The inner iteration \eqref{eq:tseng} can be written then as:
\begin{subequations}\label{eq:inner_iter}\begin{align}
   &\vec x^{(m)} = \vec u^{(m)} + \frac{\gamma}{\tau}(\vec u^k-\vec u^{(m)}) - \frac{\gamma}{\epsilon} \nabla W(\vec u^{(m)})\\
   &\vec y^{(m)} =  (id-\epsilon\gamma\Delta)^{-1} \vec x^{(m)}\\
   &\vec z^{(m)} =  \vec y^{(m)} + \frac{\gamma}{\tau}(\vec u^k-\vec y^{(m)})  - \frac{\gamma}{\epsilon} \nabla W(\vec y^{(m)})\\
   &\vec u^{(m+1)} = \operatorname{proj}_C(\vec u^{(m)} - \vec x^{(m)} + \vec z^{(m)})
\end{align}\end{subequations}
Identifying the second step of the iteration as a \emph{diffusion} step with diffusion constant $\epsilon\gamma$, 
and the first/third steps as \emph{reaction} steps driven by a chemical potential $W$,
we can conceptualize the iteration \eqref{eq:inner_iter} as a \emph{Reaction-Diffusion-Projection} scheme.
Furthermore, the second `reaction' step and the $\vec u^{(m)} - \vec x^{(m)} + \vec z^{(m)}$ term in the projection give it a \emph{Predictor-Corrector} flavor.
Given the prevalence of reaction-diffusion equations in the study of pattern formation, pioneered by A.M.~Turing in \cite{turing1990chemical}, 
this connection might not come as a complete surprise.

\section{Numerical Results}
\label{sec:num_results}

\begin{figure}
  \centering
 \includegraphics[width=0.9\textwidth]{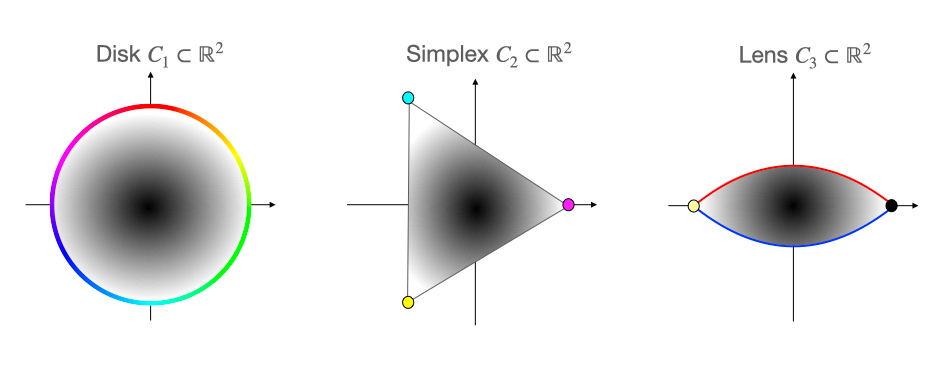}
 \caption{\emph{Obstacle sets.} The unit disk $C_1$ (left) features a continuous set of minima of the concave potential $W$ on its boundary; we use a rainbow color map to visualize them.
 The two-dimensional simplex $C_2$ (middle) has three minimal points for $W$ located at its three vertices, colored cyan/yellow/magenta respectively.
 The lens shape $C_3$ (right) has two minimal points, colored yellow/black, connected via two distinct paths (one `red', one `blue') over the boundary $\partial C_3$.}
 \label{fig:obstacles}
 \end{figure}

In this section we present a number of numerical results, computed based on the scheme \eqref{eq:inner_iter} as implemented in Algorithm \ref{alg:RDP}.
We run our calculations on the rectangular domain $\Omega=[-1,1]^2 \subset\mathbb{R}^2$, discretized via a uniform grid of 512x512 cells.
The initial conditions for each run are drawn independently for each grid cell from the uniform random distribution $U(-1,1)^2$, 
i.e.~we start each run with a maximally disordered initial state. 
We assume periodic boundary conditions which, together with the shape of the domain, 
allow us to efficiently invert the linear operator $id-\epsilon\gamma \Delta$ with the use of the Fast Fourier Transform (FFT).
Unless stated otherwise, we use the repulsive quadratic potential $W(\vec u)=\frac{1}{2\omega}(1-\lvert u\rvert^2)$,
together with one of the two-dimensional convex obstacle sets $C_1,C_2,C_3$ shown in Fig.~\ref{fig:obstacles}. 
The constants take the values $\gamma=h/20$, $\tau=h/10$, $\epsilon=3h$ (where $h$ is the grid step) and $\omega=.5$, and the tolerance is $tol=1e-6$.

The code was written in Python with the use of the libraries Numpy, Scipy (for the FFT) and Numba (for JIT optimization, especially of the projection operations).
The performance is of the order of 1 time step/sec on an Apple M1 Pro CPU.

\begin{algorithm}
  \caption{Reaction-Diffusion-Projection (single time step)}
  \label{alg:RDP}
  \DontPrintSemicolon
  \SetKwInput{KwData}{Given}
  \SetKwRepeat{Repeat}{Repeat}{until}
  \SetKw{Return}{Return}
  \BlankLine
  \KwIn{Phase field data $\vec u^k$ at time $t^k$.}
  \KwOut{Phase field data $\vec u^{k+1}$ at time $t^{k+1}=t^k+\tau$.}
  \BlankLine
  \KwData{Consts $\gamma,\tau,\epsilon,\omega$ with $0<\gamma<\tau<\epsilon\omega$, and tolerance $tol$.\\ 
  \Indp Projection operator $\operatorname{proj}_C$ onto the convex set $C$.\\
  Solver for inverting $id-\epsilon\gamma\Delta$ in domain $\Omega$ with given BCs.}
  \BlankLine
  Initialize $m=0$, $\vec u^{(0)} := \vec u^k$.\;
  \Repeat{$\lVert \vec u^{(m)} - \vec u^{(m-1)}\rVert<tol$.}{
      \linespread{1.35}\selectfont
      $\vec x^{(m)} := \vec u^{(m)} + \frac{\gamma}{\tau}(\vec u^k-\vec u^{(m)})  - \frac{\gamma}{\epsilon} \nabla W(\vec u^{(m)}) = (1-\frac{\gamma}{\tau}+\frac{\gamma}{\epsilon\omega})\vec u^{(m)} + \frac{\gamma}{\tau} \vec u^k$,\;
      solve $\vec y^{(m)} =  (id-\epsilon\gamma\Delta)^{-1} \vec x^{(m)}$,\;
      $\vec z^{(m)} := (1-\frac{\gamma}{\tau}+\frac{\gamma}{\epsilon\omega})\vec y^{(m)} + \frac{\gamma}{\tau} \vec u^k$,\;
      $\vec u^{(m+1)} := \operatorname{proj}_C(\vec u^{(m)} - \vec x^{(m)} + \vec z^{(m)})$,\;
      $m\leftarrow m+1$,\;
    }
  \Return{$\vec u^{k+1} := \vec u^{(m)}$.}
  \BlankLine
\end{algorithm}


\paragraph{Multiple Phases} 

\begin{figure}
  \centering
 \includegraphics[width=0.4\textwidth]{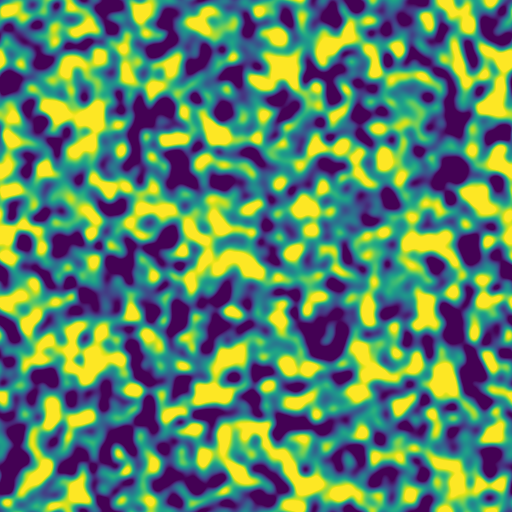}\hspace{5pt}
 \includegraphics[width=0.4\textwidth]{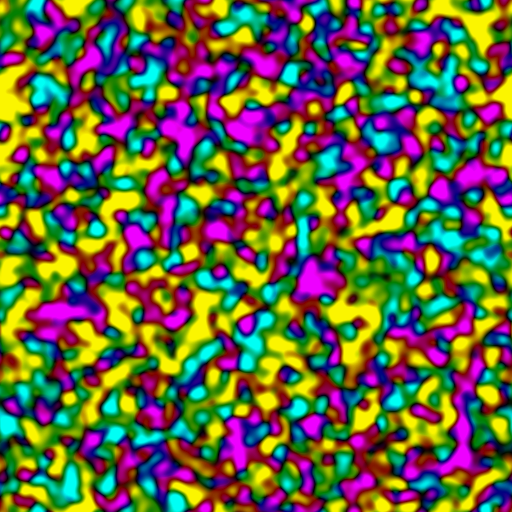}\\
 \includegraphics[width=0.4\textwidth]{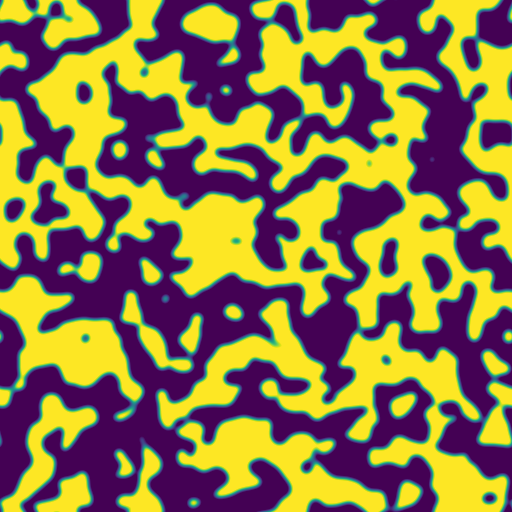}\hspace{5pt}
 \includegraphics[width=0.4\textwidth]{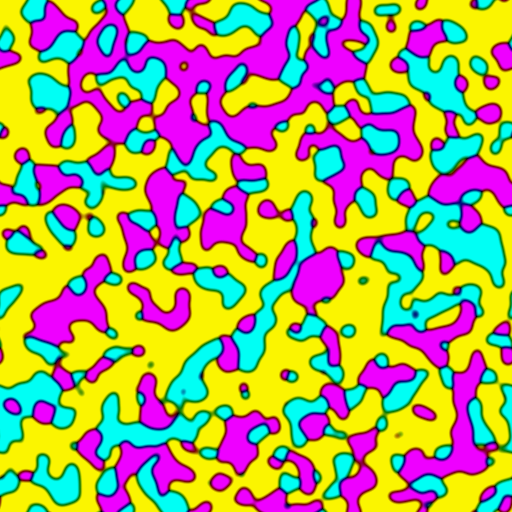}\\
 \includegraphics[width=0.4\textwidth]{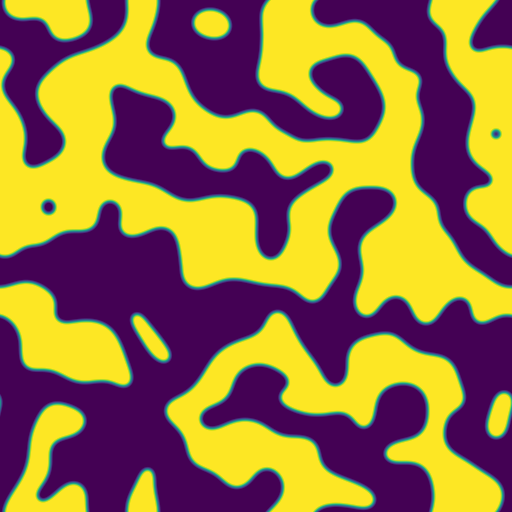}\hspace{5pt}
 \includegraphics[width=0.4\textwidth]{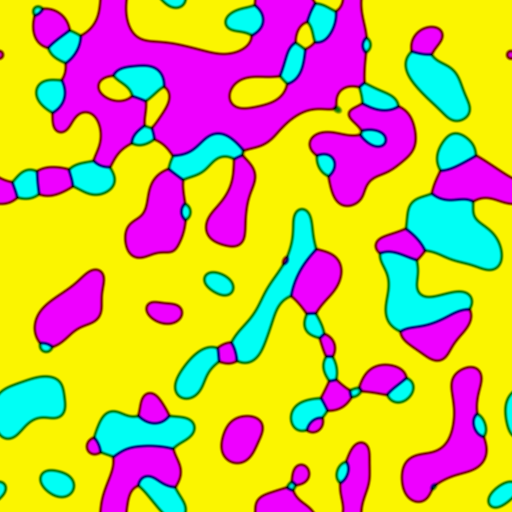}\\
 \caption{\emph{Multiple Phases.} Time evolution of 1D-valued phase field with two phases (left column) 
 and 2D-valued phase field with three phases (right column).
 Both phase fields exhibit a) separation into phases at early times, 
 followed by b) coarsening where small isolated `droplets' shrink and disappear, and finally
 c) motion of the interface by mean curvature for larger times. 
 Note the $120^\circ$ triple points in the three-phase case.}
 \label{fig:twophases}
 \end{figure}

The most straightforward application of the scheme is to reproduce the two-phases setting of the classic Allen-Cahn (Fig.~\ref{fig:twophases}, left).
The results are qualitatively similar, with the phase field quickly separating into two distinct phases,
followed by the rapid `evaporation' of the small connected components of each phase,
until the asymptotic regime of interfacial motion by mean curvature settles in. 
Compare the patterns formed at the early stages of Fig.~\ref{fig:twophases} with the \emph{dual-phase steel alloy grains} in Fig.~\ref{fig:examples},
that inspired the original Allen-Cahn model in the first place.

The main difference is that in this case, \emph{the interface has compact support}.
Indeed, going through the calculations \eqref{eq:optE}-\eqref{eq:optLambda} for $W(u)=\frac{1}{2}(1-u^2)$,
we find that the interface has sinusoidal profile of the form $U(\xi)=\sin(\xi),\xi\in[-\pi/2,\pi/2],$ 
and therefore exact thickness $\Lambda=\pi\epsilon$. 
Compare to the classic Allen-Cahn, where the interface has a diffuse $\tanh$-profile 
with approximate thickness of order $\epsilon$ and long asymptotic tails extending well into the two phases. It is noteworthy that the compact support of the interface carries over to the numerics;
since each iteration of the splitting scheme ends on a projection step, 
the vast majority of the cells away from the interface have their values projected exactly to one of the phases.

To model more than two phases interacting, we need to go to higher dimensional constraint sets. 
In the right column of Fig.~\ref{fig:twophases}, 
we show the time evolution of a phase field with values constrained in the two-dimensional simplex $C_2$ of Fig.~\ref{fig:obstacles}.
The phase field evolves through the same stages described above for the two-phases case, 
separation of phases followed by coarsening and eventually mean curvature motion of the interfaces.
One important `sanity check' for the three-phase case, which the numerical solution passes, 
is that the triple points, where all three phases meet, maintain $120^\circ$ angles while the interface moves.

In case one is wondering why we need to go to higher dimensional simplices (such as a tetrahedron in $\mathbb{R}^3$) 
to model more phases, instead of simply using a regular polygon in $\mathbb{R}^2$,
the answer is essentially a matter of topology. 
For the interface between two phases to be stable, 
there needs to exist a corresponding edge at the boundary $\partial C$ of the constraint set that connects them.
For $N$ phases to all connect to each other with edges over the boundary, 
it is easy to show that one needs to place them on an $(N-1)$--dimensional simplex.

\paragraph{Vortices} 

\begin{figure}
  \centering
 \includegraphics[width=0.8\textwidth]{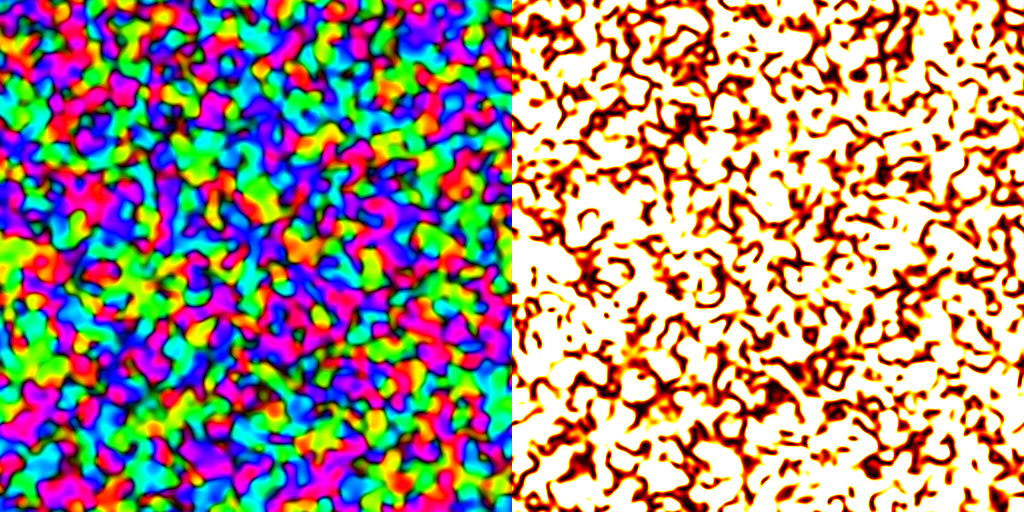}\\
 \includegraphics[width=0.8\textwidth]{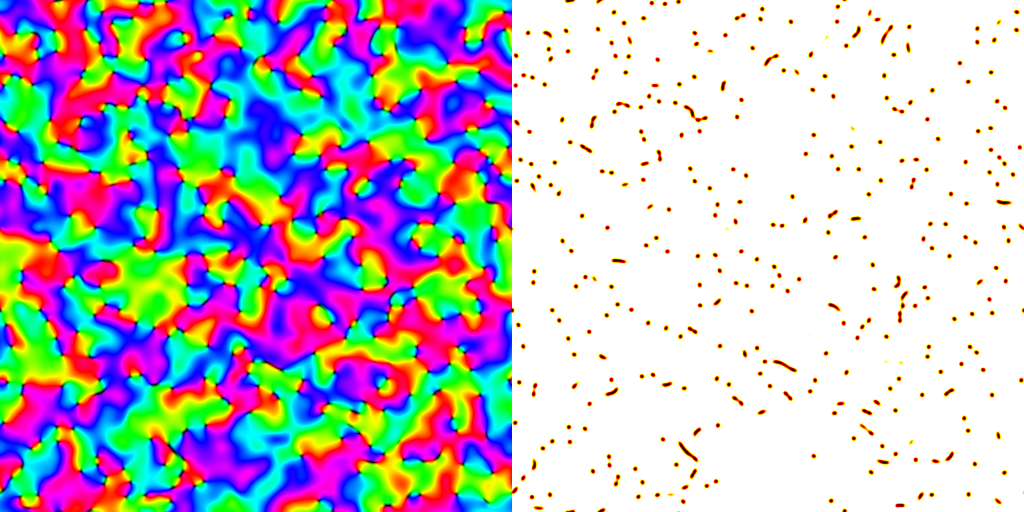}\\
 \includegraphics[width=0.8\textwidth]{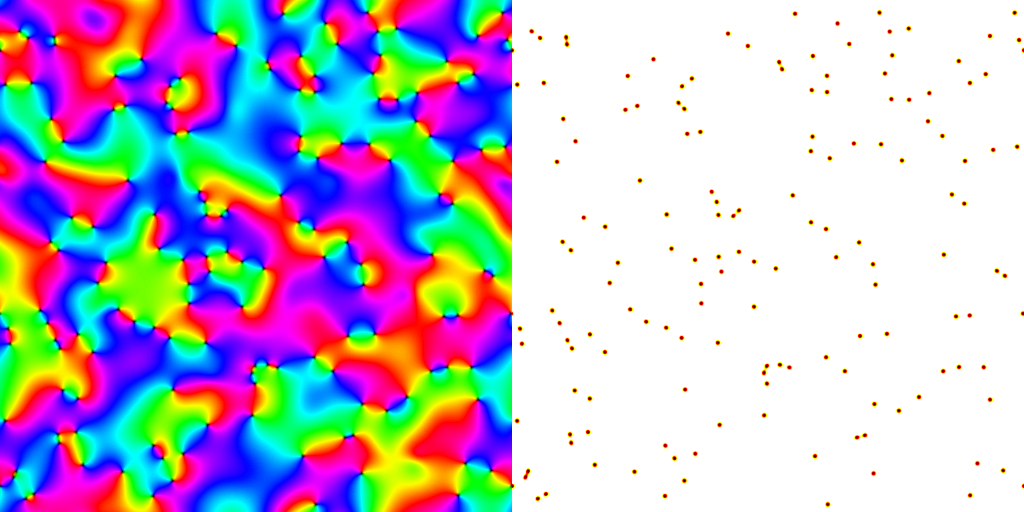}\\
 \caption{\emph{Vortices.} Time evolution of a phase field under the disk obstacle set ($C_1$ in Fig.~\ref{fig:obstacles}); 
 phase visualized as color on the left, points where $\lvert \vec u\rvert\approx 0$ localized as dark spots on the right.
 The phase field quickly converges to a harmonic map between numerous poles. 
 Pairs of nearby poles of opposite polarity converge, to their eventual mutual annihilation and the great reduction of the number of poles at later times.}
 \label{fig:vortices} 
\end{figure}

\begin{figure}
  \centering
 \includegraphics[width=0.8\textwidth]{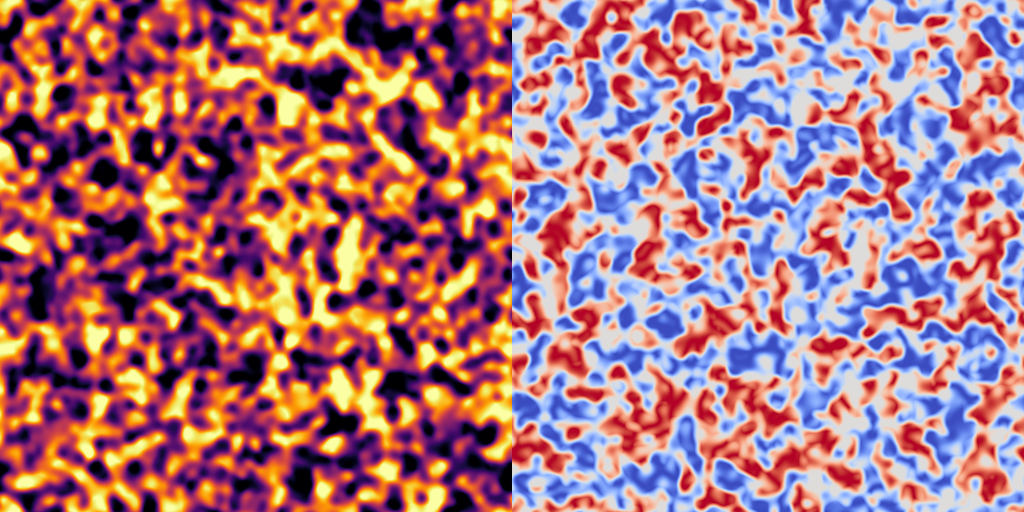}\\
 \includegraphics[width=0.8\textwidth]{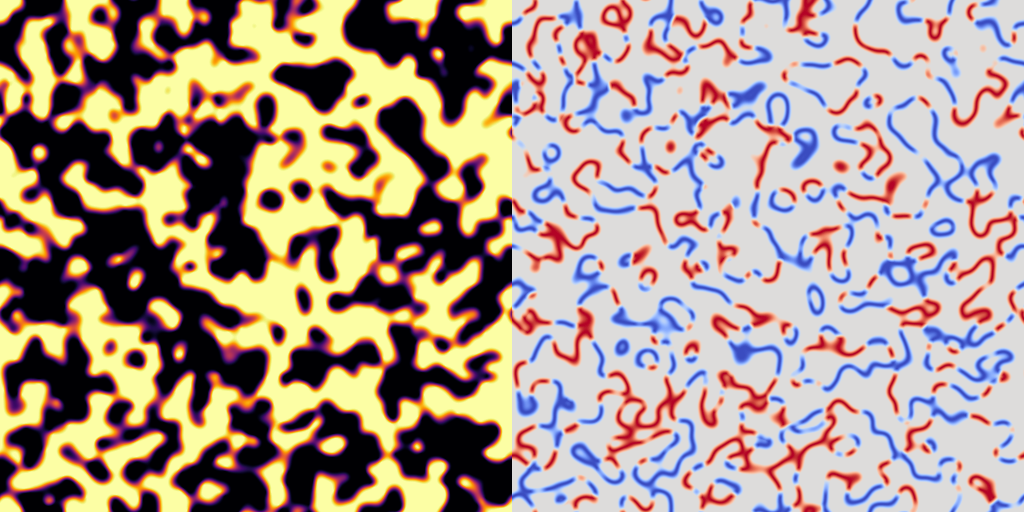}\\
 \includegraphics[width=0.8\textwidth]{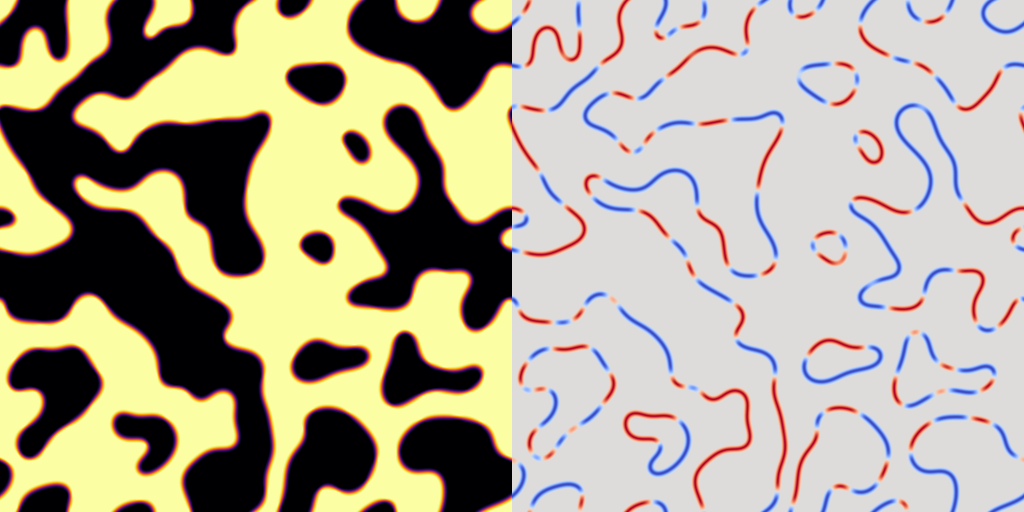}\\
 \caption{\emph{Nested Interfaces.} Time evolution of a phase field under the lens obstacle set ($C_3$ in Fig.~\ref{fig:obstacles}). 
 The phase field quickly settles into a mix of two phases, visualized as yellow/black on the left, similar to the two-phases case of Fig.~\ref{fig:twophases}.
 The interface between the two phases (visualized on the right) is itself partitioned into two phases, 
 corresponding to the two extremal paths (upper/red, lower/blue) shown in Fig.~\ref{fig:obstacles}.
 }
 \label{fig:nested}
 \end{figure}

In Fig.~\ref{fig:vortices}, the evolution of a phase field with values constrained inside the unit disk ($C_1$ in Fig.~\ref{fig:obstacles}) is presented.
Contrasted against the results discussed in the previous section, we can immediately see the difference between 
isolated and connected sets of minimal points on $\partial C$.
Continuous sets of minimal points act as a single non-uniform phase;
within each patch of this phase, and as discussed in section\ref{sec:theory}, the phase field tries to `relax' into a smooth harmonic map.
Whenever for topological reasons this is impossible globally, 
such as for a map from a periodic rectangular domain onto the non-simply connected unit circle in this case,
the phase field is `frustrated' and reacts by localizing the topological defects as much as possible.
In Fig.~\ref{fig:vortices}, this takes the form of point-like vortices of two possible polarities,
as the phase field can map onto the unit circle clockwise or counterclockwise in their vicinity.
These vortices exhibit emergent long-term dynamics, where neighboring vortices of opposite polarity are 
attracted to and ultimately cancel each other out.
In this manner the phase field is slowly `untangling' itself, reducing the total energy in distinct jumps with each annihilation event.
This version of the model is closely related to the behavior of \emph{thin liquid crystal films}, as seen in Fig.~\ref{fig:examples}.

\paragraph{Nested Interfaces} 

The dynamics of Fig.~\ref{fig:vortices}, 
where the phase field is constrained to take values from the lenticular set $C_3$ of Fig.~\ref{fig:obstacles},
serve as an even stronger demonstration of the influence of the shape of $C$ on the observed behavior.
It is a return to the two isolated minima of the very first case discussed, 
but this time there are two separate edges connecting the two minima over the boundary of $C$,
shown as red/upper and blue/lower in Fig.~\ref{fig:obstacles}.
As a consequence, the interface between bulk phases (drawn in yellow and black in Fig.~\ref{fig:vortices}) can itself exist in two distinct stable states.
We observe then two simultaneous effects: the motion by mean curvature of the interfaces,
coupled with the coarsening of the interfacial phases within each interface.

At this point we are reaching the limits of what can be observed over a two-dimensional domain;
where we to repeat this calculation over a three-dimensional domain,
the interfaces between bulk phases would be curved surfaces,
themselves partitioned into evolving red/blue patches.
The interfaces between these patches would be curves following motion by curvature, 
embedded in surfaces that are also evolving under mean curvature,
corresponding to an interesting set of coupled geometric PDEs.
This configuration has also important connections to the study of \emph{biological cell membranes},
since the bilipid layer that forms these membranes can exist in two different states \cite{sych2021does} 
(see Fig.~1 for a relevant image from this paper).
Further analysis and numerical exploration of this interesting case are left as potential future work.

\paragraph{Periodic Structures} 

\begin{figure}
  \centering
 \includegraphics[width=0.4\textwidth]{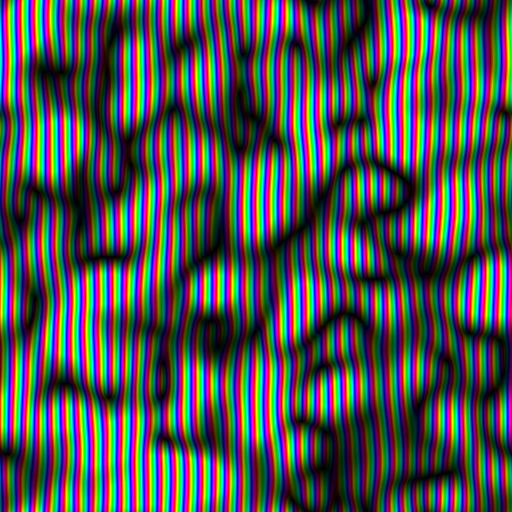}
 \includegraphics[width=0.4\textwidth]{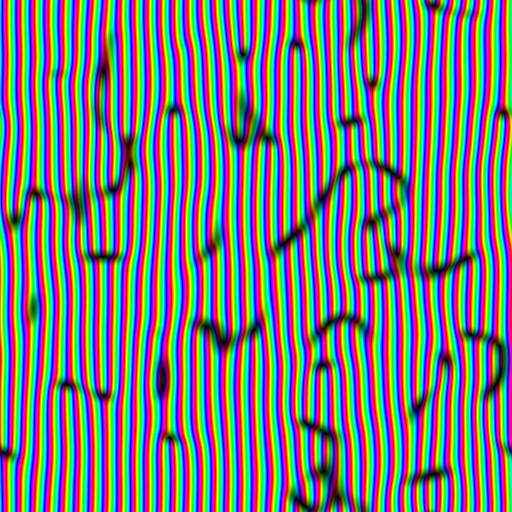}\\[2pt]
 \includegraphics[width=0.4\textwidth]{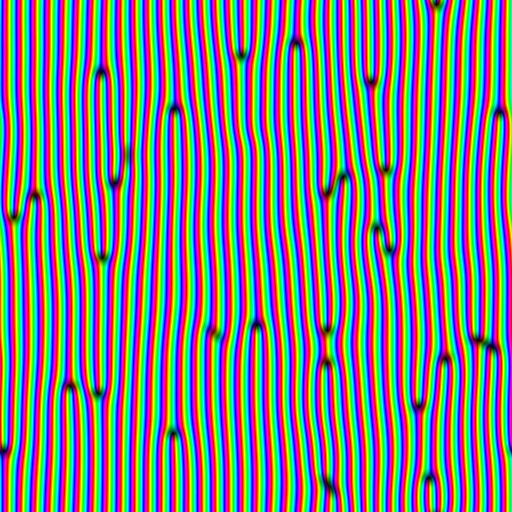}
 \includegraphics[width=0.4\textwidth]{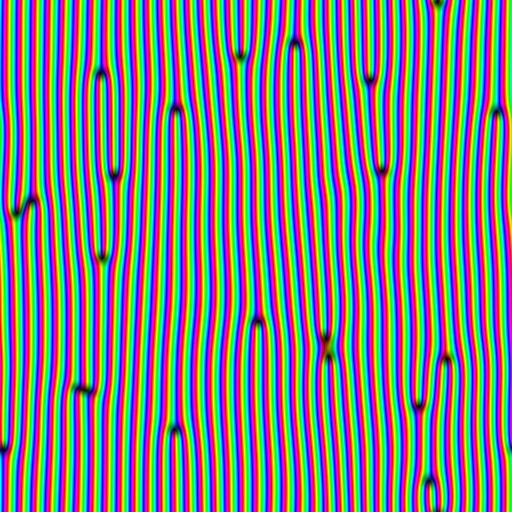}\\[2pt]
 \caption{\emph{Periodic Structures.} Phase field evolving under the disk obstacle set ($C_1$ in Fig.~\ref{fig:obstacles}), 
 but with an anisotropic variant of the gradient term $\lvert\nabla\vec u\rvert^2$ in the energy which promotes periodicity in the $x$-direction.
 As the field settles into periodic patches that are out-of-phase with each other, 
 the topological defects in the pattern consolidate into points. 
 These `kinks' slowly converge and annihilate each other similar to the poles in Fig.~\ref{fig:vortices}.}
 \label{fig:periodic}
 \end{figure}

A useful observation about the generalized Ginzburg-Landau functional \eqref{eq:GL} is that it is \emph{isotropic},
as it is invariant under rotations $\vec{u}\rightarrow R_\theta \vec{u}$.
In Fig.~\ref{fig:periodic}, we present the behavior of an \emph{anisotropic} variant of the functional 
with the modified smoothness term:
\begin{equation}
 \lvert\nabla u_1\rvert^2 + \lvert\nabla u_2\rvert^2 
 = \left\lvert\frac{\partial\vec u}{\partial x}\right\rvert^2 + \left\lvert\frac{\partial\vec u}{\partial y}\right\rvert^2
 \quad\longrightarrow\quad \left\lvert\frac{\partial\vec u}{\partial x} - R_{90^\circ}\vec u\right\rvert^2 + \left\lvert\frac{\partial\vec u}{\partial y}\right\rvert^2 \,,
\end{equation}
where $R_{90^\circ} = \left(\begin{smallmatrix} 0 & -1\\1 & 0\end{smallmatrix}\right)$ is the $90^\circ$-counterclockwise rotation matrix.
The anisotropic term promotes uniform values in the $y$-direction, same as the isotropic one, 
but it favors periodic values in the $x$-direction.
Mutatis mutandis, in particular by replacing the operator $id-\epsilon\gamma\Delta$:
\begin{equation}
  id-\epsilon\gamma\Delta \quad\longrightarrow\quad (1+\epsilon\gamma)id - \epsilon\gamma\Delta + 2R_{90^\circ}\partial_x\,,
 \end{equation}
the Reaction-Diffusion-Projection scheme (Alg.~\ref{alg:RDP}) works as presented.
We constrain the phase field to take values in the unit disk ($C_1$ in Fig.~\ref{fig:obstacles}), as in the vortices case, 
and observe, as expected, an initial synchronization of the phases into local aligned patches separated by `cracks'.
The alignment phase continues until the topological defects in the pattern focus into points.
In the long run, these behave similar to the vortices of Fig.~\ref{fig:vortices}, 
in that they come in two orientations, `(forking) up' and `(forking) down', 
that are attracted to and eventually cancel each other.
This type of periodic linear pattern is fairly common in nature, perhaps most recognizably in \emph{sand waves} (Fig.~\ref{fig:examples}).

 \paragraph{Traveling Waves} 

\begin{figure}
  \centering
 \includegraphics[width=0.4\textwidth]{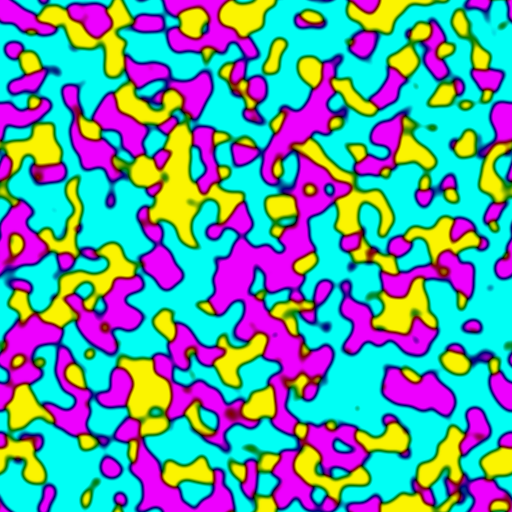}
 \includegraphics[width=0.4\textwidth]{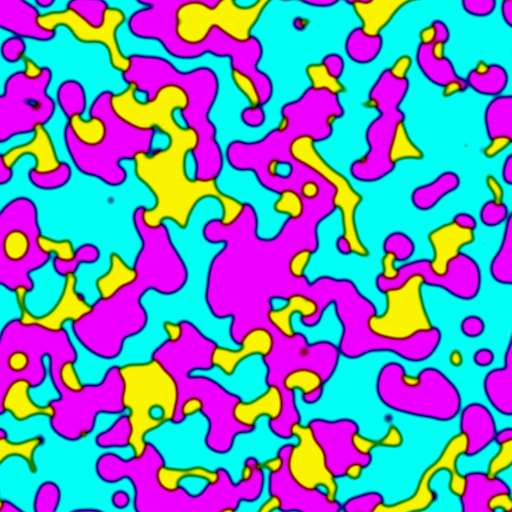}\\[2pt]
 \includegraphics[width=0.4\textwidth]{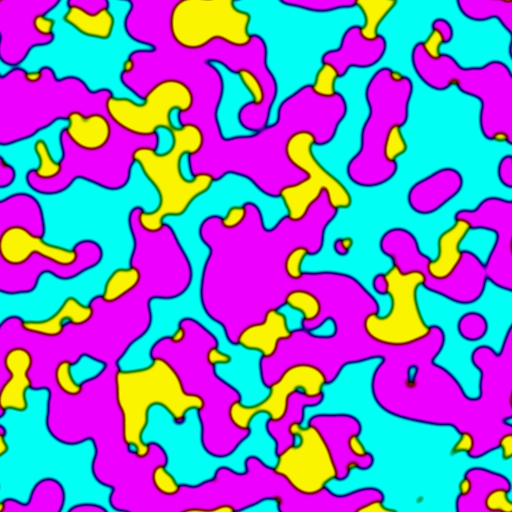}
 \includegraphics[width=0.4\textwidth]{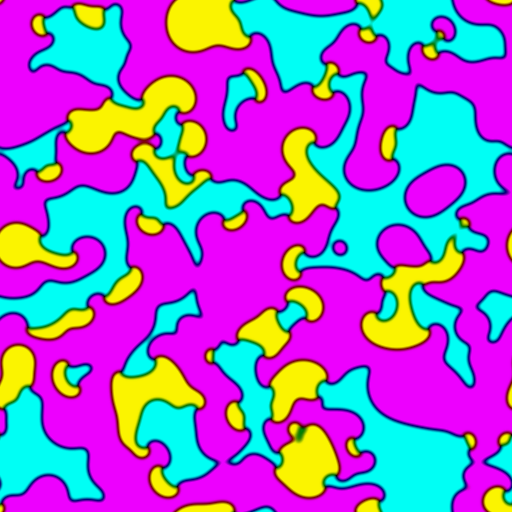}\\[2pt]
 \caption{\emph{Traveling Waves.} Phase field evolving under the simplex obstacle set ($C_2$ in Fig.~\ref{fig:obstacles}), 
 but with a non-conservative variant of the driving force $\nabla W(\vec u)$. 
 The variant force makes each of the three phases unstable in the presence of one of the others in a rock-paper-scissors pattern;
 magenta `beats' cyan `beats' yellow `beats' magenta.
 Hence, each interface where two phases meet tends to advance into the `losing' phase of the pair in a traveling wave fashion.}
 \label{fig:twave}
 \end{figure}

 In Fig.~\ref{fig:twave} we consider another variation of the Reaction-Diffusion-Projection scheme (Alg.~\ref{alg:RDP}).
 We use the triangle obstacle set $C_1$ of Fig.~\ref{fig:obstacles}, as in the three-phases case, 
 but we replace the driving force $\nabla W(\vec u)$ of the reaction step(s) with a non-conservative one, 
 i.e.~one that is not a gradient of a scalar potential:
\begin{equation}
  -\gamma\nabla W(\vec u) = \frac{\gamma}{\omega}\vec u \quad\longrightarrow\quad \frac{\gamma}{\omega}R_{30^\circ}\vec u
\end{equation}
The splitting algorithm remains valid under this alteration, 
as Tseng's scheme \eqref{eq:tseng}, on which it is built, 
is in fact a particular case of a more general scheme (Prop.~25.36 in \cite{bauschke2011convex}).
That version of the scheme effectively allows us to replace the conservative term $\mathcal{G}'(\vec u)$ 
with a monotone operator $B$, 
such that $(\vec u-\vec v)\cdot(B\vec u-B\vec v)\geq 0$ for any $\vec u$,$\vec v$. 
The rotation $R_{30^\circ}$, in fact any rotation $R_\theta$ with $\theta\leq 90^\circ$, is indeed such a monotone operator.

Referring again to the triangular set $C_2$ in Fig.~\ref{fig:obstacles}, 
the conservative force $\nabla W$ points radially outwards, 
so that any point on the edges of the set is forced to 'slide' towards the vertex closest to it.
The vertices themselves are stable, resulting in the three-phase case of Fig.~\ref{fig:twophases}.
On the other hand, the introduction of the 30$^\circ$ counterclockwise rotation is chosen 
so that at each vertex of the triangle the force is perpendicular to one of the adjacent edges.
This makes the vertex metastable, as with even the slightest perturbation in the counterclockwise direction 
the rotated force starts pushing towards the next vertex. 
The results can be seen in Fig.~\ref{fig:twave}; the three phases destabilize each other in a rock-paper-scissors fashion,
magenta `beats' cyan `beats' yellow `beats' magenta. 
Wherever two phases meet, the interface advances into the `losing' phase in the form of a traveling wave
reminiscent of the Belousov-Zhabotinsky chemical reaction (see Fig.~\ref{fig:examples}).

 \section{Future Directions}
 \label{sec:future}

 As mentioned at various points above, there are a number of exciting possibilities for further work.
 One straightforward direction is to \emph{couple with other equations} from Physics, Chemistry and Biology,
 as has been done very fruitfully for the classic phase field method. 
 The fact that, with the scheme presented here, the interfaces are compact and the phases take exact values,
 could positively impact the behavior of the coupled systems. 
 Furthermore, one can use the demonstrated ability of the scheme to reproduce features of physical systems,
 that would otherwise need to be captured via for instance additional reaction-diffusion equations, 
 to derive more compact mathematical models.

 Significant work can also be done in the direction of \emph{improving the performance} of the scheme in practice,
 especially in the context of doing numerics on three-dimensional domains. 
 The compactness of the interface, especially at later times when the coarsening has progressed sufficiently,
 lends itself well to the \emph{narrow band method}, 
 popularized in the level sets context \cite{adalsteinsson1995fast} but applicable to PDEs in general \cite{deckelnick2010h}.

 A more theoretical approach would be \emph{to study the effect of the various possible constraint sets} more systematically.
 The exact nature and dynamics of the various interfaces and topological singularities,
 the limit to the associated geometric PDEs, and the asymptotic behavior in a multitude of time and space scales,
 have all proven very fruitful research areas in the classic phase field context.
 Powerful mathematical tools, such as $\Gamma$-convergence and formal asymptotics, 
 were used successfully there, and could be deployed here too with potentially very interesting results.

 A final intriguing idea is to \emph{apply the scheme to non-Euclidean domains}, 
 such as curved surfaces (discretized as triangle or quad meshes), or even completely general graphs.
 A quick survey of the Reaction-Diffusion-Projection scheme (Alg.~\ref{alg:RDP}) 
 reveals that the Reaction and Projection steps are completely local, and can be performed node-wise.
 The Diffusion step on the other hand is non-local, 
 but its calculation has attracted extensive attention 
 on discrete surfaces from the Computational Geometry and Graphics communities \cite{vaxman2010multi} \cite{crane2013geodesics},
 and more recently on graphs by the Machine Learning community \cite{gasteiger2019diffusion}.
 The computational tools to make the transition are therefore widely available,
 and the results could be presented to these very active communities.




  \bibliographystyle{elsarticle-num-names} 
  \bibliography{vantzos2023}


%
%
%
\end{document}